\documentclass[letterpaper, paper,11pt]{AAS}	

\usepackage{bm}
\usepackage{amsmath}
\usepackage{amssymb}
\usepackage{subfigure}
\usepackage{subcaption}
\usepackage[colorlinks=true, pdfstartview=FitV, linkcolor=black, citecolor= black, urlcolor= black]{hyperref}
\usepackage{overcite}
\usepackage{footnpag}			      	
\usepackage{graphicx}%
\usepackage{multirow}%
\usepackage{amsmath,amssymb,amsfonts}%
\usepackage{amsthm}%
\usepackage{mathrsfs}%
\usepackage[title]{appendix}%
\usepackage{xcolor}%
\usepackage{textcomp}%
\usepackage{manyfoot}%
\usepackage{booktabs}%
\usepackage{algorithm}%
\usepackage{algorithmicx}%
\usepackage{algpseudocode}%
\usepackage{listings}%
\usepackage{multirow}

\usepackage{bm}
\usepackage{gensymb}
\usepackage{hyperref}
\usepackage{nomencl}
\usepackage{caption}
\usepackage{subcaption}
\usepackage{tabularx}
\usepackage{multirow}
\usepackage{xcolor}
\usepackage{amsmath}
\hypersetup{hidelinks} 
\usepackage{algorithm}
\usepackage{algpseudocode}
\usepackage{multicol}

\PaperNumber{24-426}

\begin{document}

\title{Time-triggered Reduced Desensitization Formulation For Solving Optimal Control Problems
}

\author{Praveen Jawaharlal Ayyanathan\thanks{PhD Candidate, Department of Aerospace Engineering, Auburn University, Auburn, AL 36849, USA.},  
Ehsan Taheri\thanks{Assistant Professor, Department of Aerospace Engineering, Auburn University, Auburn, AL 36849, USA.}
}

\maketitle{}

\begin{abstract}
Fuel-optimal trajectories are inherently sensitive to variations in model parameters, such as propulsion system thrust magnitude. This inherent sensitivity can lead to dispersions in cost-functional values, when model parameters have uncertainties. Desensitized optimal control aims at generating robust optimal solutions while taking into account uncertainties in the model parameters. While desensitization techniques typically apply along the entire flight time, this paper introduces a novel time-triggered desensitization mechanism by modifying a recently developed desensitization method -- the Reduced Desensitization Formulation (RDF). By selectively desensitizing over specific time intervals of trajectories, we demonstrate the improved optimality of desensitized trajectories.
We investigate the effects of temporal desensitization on the final cost and trajectory by considering thrust magnitude uncertainty for two classes of low-thrust trajectory optimization problems: 1) minimum-fuel rendezvous maneuvers and 2) orbit-raising maneuvers. Results show that temporal desensitization can achieve similar dispersion levels to full mission desensitization with an improved final cost functional. 
\end{abstract}

\section{Introduction}
Optimal control problems (OCPs) are subject to uncertainties in system parameters, and will inevitably, exhibit dispersions in solutions in the presence of model uncertainties. The inherent sensitivity of the optimal solutions to model uncertainties has a practical significance, since many system parameters can change due to various factors in real-world applications. From a practical point of view, uncertainty-aware trajectory optimization can be useful for efficient and robust operation of unmanned aerial vehicles and space vehicles, since they operate within unpredictable environments with system parameter variations. Desensitized Optimal Control (DOC) \cite{seywald2003desensitized} belongs to the class of robust control methods that seeks to reduce the sensitivities of the set of solutions to optimal control problems (OCPs) with respect to the uncertain parameters of the dynamical model. 

Several techniques exist in the literature that formulate and achieve desensitization through different approaches. For instance, in \cite{seywald2024desensitized}, the authors represent the sensitivity of the physical states at the final time with respect to perturbations in the physical states at the current time through a Lagrange multiplier-like quantity to achieve desensitization.  A comparison between different sensitivity calculation methods, that can be used for the desensitization of OCPs, is presented in \cite{akman2020efficient}. Reference \cite{makkapati2020c}, introduces a costate-based desensitization technique by considering the costates of uncertain parameters. The uncertain parameter's costate is then penalized along the trajectory to achieve desensitization. Penalization is typically achieved by introducing a quadratic penalty term to the running cost of standard (i.e., sensitive) OCPs. 

In our recent work \cite{jawaharlal2024reduced}, we developed a computationally efficient desensitization formulation, called the Reduced Desensitization Formulation (RDF) that leverages the costate-desensitization idea, but significantly reduces the number of differential equations required for desensitizing the solutions. Further, the utility of the RDF is shown by applying it to three classes of problems: 1) fuel-optimal interplanetary trajectory optimization problems, 2) fuel-optimal rocket-landing trajectory optimization problems, and 3) maximum-final-radius orbit-raising problems. Other approaches to robust trajectory design include polynomial chaos expansion based trajectory optimization for uncertainty analysis \cite{li2014aircraft}. Reference \cite{wang2022desensitized} introduces a DOC strategy for electric vehicles to reduce the sensitivity of battery pack temperature. The advantages of sensitivity (both forward sensitivity and adjoint method) for robust re-entry trajectory optimization for an Apollo-like capsule are compared in Ref.\cite{akman2020efficient}.  

Trajectory desensitization offers a potential approach to lower mission operational expenses by minimizing the need for frequent trajectory re-planning. This is achieved by designing trajectories that exhibit reduced sensitivity to uncertain parameters. The importance of considering operational aspects in spacecraft guidance, navigation, and control design is emphasized in reference \cite{morris2011benefits}. In the case of satellite formation flying missions, thrust uncertainties have been identified as a significant challenge\cite{zhang2014nanosatellite}. Many researchers have addressed the uncertainties that are encountered during space missions. In Ref.~\cite{burnettdesensitized}, the authors have developed desensitized robust guidance techniques for spacecraft formation flying around asteroids, accounting for uncertainties in dominant gravitational harmonics and solar radiation pressure disturbances. Furthermore, probability-based trajectory optimization methods have been explored for close proximity operations, taking into account uncertainties in path constraints \cite{jewison2018probabilistic}. 

DOC strategies are also considered for Entry, Descent and Landing (EDL) problems.  The importance of pinpoint landing on Mars has been a focus of numerous studies \cite{li2014review, wan2022fuel} leading to the development of various robust trajectory generation methods. For instance, in \cite{xu2024reduced}, the authors propose a reduced-order desensitization method for generating robust reference trajectories for the Entry Terminal Point Control (ETPC) algorithm. Probability-based approaches have also been explored for Mars atmospheric entry, with an uncertainty considered in the initial states and atmospheric density \cite{yu2016observability}. In \cite{akman2018using}, the authors desensitize an optimal re-entry trajectory with respect to air density fluctuations. 

In Reference \cite{shen2010desensitizing}, the sensitivity-based DOC strategy has been used to minimize fuel consumption and reduce the sensitivity to state uncertainties for a powered descent Mars pinpoint landing problem. Direct collocation and nonlinear programming techniques are also used to solve Mars entry trajectory optimization problems using the sensitivity-based DOC method \cite{li2011mars}. It is shown, that a reduction in sensitivity of terminal state variables with respect to uncertainties on initial state variables and Mars atmospheric density variations is achievable. The sensitivity-matrix-based DOC method is also applied to landing problems with uncertainties in atmospheric density and aerodynamic characteristics \cite{xu2015robust}. The landing error on small bodies is reduced in the presence of uncertainties of target body and thrust error in \cite{hu2016desensitized} using a desensitization formulation.

Although the idea of desensitization is advantageous in reducing the sensitivity of the optimal solution to uncertainties, it has an adverse impact on the value of the cost (performance index) of the original, sensitive OCP. In other words, there is an inherent trade-off between optimality and robustness \cite{jawaharlal2024reduced} and a balance between the level of desensitization and optimality is an important trade-off that has to be taken into account for real-world applications of the DOC methods. 
\newpage
The main contribution of this paper is to develop a methodology to provide answers to the following questions: 
\begin{itemize}
    \item 1) Is it possible to improve the final cost (performance index) of a desensitized trajectory?
    \item 2) Should OCPs be desensitized along their entire trajectory? or Is it beneficial to isolate certain intervals and/or regions of interest and reduce the solution sensitivity only during those intervals/regions?
    \item 3) What is the impact of time-triggered desensitization on fuel-optimal and orbit-raising trajectory optimization problems?
\end{itemize}

In this paper, we address the above-mentioned questions and build upon the framework that we presented in \cite{jawaharlal2024reduced}. More specifically, we modify the cost functional of the RDF method by introducing a time-triggered factor that allows us to consider desensitization by selectively considering only one (or multiple non-overlapping) interval(s). 

The remainder of the paper is organized in the following manner. First, the formulation for the sensitive fuel-optimal trajectory optimization problem is presented.  Then, the time-triggered fuel-optimal RDF formulation is explained. This is followed by the time-triggered orbit-raising RDF formulation. Then, the results obtained using the proposed time-triggered RDF formulations are presented. Finally, concluding remarks are given. 

\section{Sensitive fuel-optimal trajectory optimization formulation} \label{sec:sensitive}
Consider the task of designing fuel-optimal space trajectories. The cost functional for a fixed-time, fuel-optimal trajectory optimization problem can be written as,
\begin{equation} \label{eq:costfuncTPBVP}
     \underset{\delta \in [0,1] \& || \hat{\bm{u}}|| = 1}{\mathop{\text{minimize}}} ~J = -m(t_f),
\end{equation}
where $\delta \in [0,1]$ is the engine throttling input, and $\hat{\bm{u}}$ is the thrust steering unit vector. Both $\delta$ and $\hat{\bm{u}}$ are control inputs. In Eq.~\eqref{eq:costfuncTPBVP}, $m(t_f)$ is the spacecraft mass at the final time. The spacecraft dynamics, using different coordinate/element sets (e.g., Cartesian coordinates or Modified Equinoctial Elements (MEEs)) \cite{taheri2016enhanced,junkins2019exploration}, can be written in a unified control-affine form as,
\begin{align} \label{eq:EOM}
    \bm{\dot{x}} & = \bm{f}(\bm{x},\bm{\Delta}_c,t) = \bm{A}(\bm{x},t) + \mathbb{B}(\bm{x},t) \bm{\Delta}_c , & \dot{m}  & = -\frac{T}{c} \delta,
\end{align}
where $\bm{x} \in \mathbb{R}^6$ denotes the coordinate/element vector, $m$ is the spacecraft mass, and $c = I_\text{sp} g_0$ denotes the effective constant exhaust velocity. Here, $I_\text{sp}$ is the specific impulse and $g_0$ is the acceleration due to gravity at sea level. In Eq.~\eqref{eq:EOM}, $T$ denotes the maximum thrust magnitude of the propulsion system. In addition, the control acceleration vector is expressed as,
\begin{align*}
    \bm{\Delta}_c = \frac{T}{m} \delta \hat{\bm{u}}.
\end{align*}

In reality and for solar-powered electric propulsion systems, the magnitude of thrust and specific impulse are functions of the power \cite{TAHERI2020166, arya2021composite}. To simplify the problem formulation, it is assumed that the values for $T$ and $c$ remain constant along the entire trajectory. In addition, perturbations due to third bodies are not considered. In this paper, the set of MEEs is used to represent the spacecraft's equations of motion (a derivation for the set of MEEs can be found in \cite{jawaharlal2023mapped}). Note that in this case, $\hat{\bm{u}}$ is the thrust steering unit vector expressed in the local-vertical local-horizontal coordinate system associated with the set of MEEs, i.e., $\hat{\bm{u}} = [u_r,u_t,u_n]^{\top}$.

Let $\bm{\lambda_x} \in \mathbb{R}^6$ denote the costate vector associated with the state vector $\bm{x}$ and let $\lambda_m$ denote the costate associated with mass. Let $H = L + \bm{\lambda_x}^{\top} \bm{f}$ denote the (optimal control) Hamiltonian. From Eq.~\eqref{eq:costfuncTPBVP}, we have $L = 0$. The costate differential equations can be obtained using the Euler-Lagrange equation as,
\begin{equation}\label{eq:costateTPBVP}
\dot{\boldsymbol{\lambda}}_{\boldsymbol{x}}=-\left[\frac{\partial H}{\partial \boldsymbol{x}}\right]^{\top}, \quad \dot{\lambda}_m=-\frac{\partial H}{\partial m}.
\end{equation}

The expressions for optimal/extremal controls can be obtained by using the Pontyagin's Minimum Principle (PMP) and Lawden's primer vector theory as \cite{jawaharlal2024reduced},
\begin{align}\label{eq:controlTPBVP}
    \hat{\bm{u}}^* & = - \frac{\mathbb{B}^{\top} \bm{\lambda}_x}{|| \mathbb{B}^{\top} \bm{\lambda}_x ||}, & \delta^*= \begin{cases}
      1, & \text{if} \ S > 0, \\
      0, & \text{if} \ S < 0, \\
\end{cases} 
\end{align}
where $S$ is the \textit{thrust switching function} and is defined as,
\begin{equation}
    S = \frac{c ||\mathbb{B}^{\top} \bm{\lambda}_x||}{m} + \lambda_m.
\end{equation}

Upon considering fixed-time, rendezvous-type low-thrust maneuvers, the initial and final states are known. Since the final value of mass is optimized, it is unknown and its costate value can be determined using the transversality condition (i.e., $\lambda_m(t_f)= -1$). The boundary conditions for the sensitive fuel-optimal trajectory optimization problem can be summarized as, 
\begin{align} \label{eq:bcTPBVP}
    \bm{\Psi}(\bm{\eta}) = \left[ \ [\bm{x}(t_f) - \bm{x}_T]^\top, \lambda_m(t_f) + 1 \right]^\top = \bm{0},
\end{align}
where $\bm{x}_T$ denote the final/target states and $\bm{\eta} = [\bm{\lambda_x}^{\top}(t_0),\lambda_m(t_0)]^{\top} \in \mathbb{R}^7$ denotes the unknown initial vector of costates. The state differential equations, given in Eq.~\eqref{eq:EOM}, costate differential equations, given in 
 Eq.~\eqref{eq:costateTPBVP}, extremal controls, given in 
 Eq.~\eqref{eq:controlTPBVP}, and the boundary conditions, given in 
 Eq.~\eqref{eq:bcTPBVP} represent the Hamiltonian two-point boundary-value problem (TPBVP) associated with the standard/sensitive OCP. The resulting TPBVPs are typically solved using single- or multiple-shooting solution schemes using quasi-Newton non-linear root-finding solvers (e.g., MATLAB's \texttt{fsolve} built-in function). 
 
 It is known that fuel-optimal trajectories consist of a finite number of throttle switches due to the optimality criterion, as it is depicted in the optimal throttle logic in Eq.~\eqref{eq:controlTPBVP}. The piece-wise continuous optimal throttle structure is a source of difficulty when numerical methods are used, since the Jacobian of the residual vector (i.e., Eq.~\eqref{eq:bcTPBVP}) with respect to the unknown decision values (i.e., $\bm{\eta}$) becomes singular \cite{haberkorn2004low}. The existence of multiple throttle switches also reduces the basin of solution of the numerical methods. For these reasons, the resulting TPBVPs are typically solved through various regularization methods, such as those that regularize the thrust/throttle non-smooth profile. Control regularization methods such as logarithmic smoothing \cite{bertrand2002new}, extended logarithmic smoothing \cite{taheri2016enhanced}, hyperbolic tangent smoothing \cite{taheri2018generic}, and L2-norm based regularization \cite{kovryzhenko2023vectorized,taheri2023l2} are among the popular and efficient methods for attenuating the issues with non-smooth control profiles. Through regularization methods, the non-smooth OCP is embedded into a one- or multiple-parameter family of smooth neighboring OCPs. Numerical continuation and/or homotopy methods are used for solving the resulting family of the OCPs until a solution to the original OCP is obtained \cite{trelat2012optimal}.  

\section{Time-Triggered fuel-optimal RDF formulation}\label{sec:time_state_rdf}
In our earlier work \cite{jawaharlal2024reduced}, we developed a desensitization methodology called the RDF using non-uniqueness of the solution of the costate differential equations. We used a hybrid indirect-direct optimization method for solving the resulting desensitized OCPs. More specifically, the indirect formalism of optimal control theory is only used for deriving the necessary differential equations of the costates. However, we do not use the extremal control expressions of the indirect method, and instead, we use a direct optimization method to solve the resulting DOC problems. 

To achieve desensitization of the cost functional with respect to thrust magnitude, thrust is elevated as a state with its time derivative set to $0$, i.e., $\Dot{T}=0$. The time rate of change of the thrust costate can be obtained using the Euler-Lagrange equation as,
\begin{align} \label{eq:thrustcostate}
   \dot{\lambda}_T = -\frac{\partial H}{\partial T}, \quad \text{with} \quad \lambda_T(t_f) = 0,
\end{align}
where $H$ is the Hamiltonian. By following the RDF in \cite{jawaharlal2024reduced}, the thrust costate differential equation is written as,
\begin{equation}\label{eq:thrustcostateMEE}
    \dot{\lambda}_{T} = -\frac{\delta}{m}(K_{vr}u_r +  K_{vt}u_t + K_{vn}u_n -  K_{m}\frac{m}{c}), 
\end{equation}
where $K_{vr}=1$, $K_{vt}=1$, $K_{vn}=1$ denote the constant velocity costates. The reason for setting all the costates to a constant value of $1$ is outlined in \cite{jawaharlal2024reduced}. As mentioned before, $u_r$, $u_t$, and $u_n$ denote radial, transversal and normal components of the thrust steering unit vector that is expressed in the osculating coordinate attached to the spacecraft. 

For desensitizing the trajectory, the thrust costate is minimized along the entire trajectory. Therefore, the cost functional in Eq.~\eqref{eq:costfuncTPBVP} is augmented with an along-the-path quadratic penalty term as,
\begin{equation} \label{eq:costfuncaug}
     \underset{\delta \in [0,1] \& || \hat{\bm{u}}|| = 1}{\mathop{\text{minimize}}} ~J = -m(t_f) + \underbrace{\int_{t_0}^{t_f}  Q \lambda^2_T(t) \ dt}_{\text{desensitization penalty term}} ,
\end{equation}
where $t_0$ is the initial time, $t_f$ is the final time and $Q$ is a positive constant parameter. When $Q=0$ the Lagrange cost (i.e., $L = Q \lambda^2_T(t)$) in Eq.~\eqref{eq:costfuncaug} becomes $0$ and the OCP corresponds to the sensitive OCP. To achieve desensitized solutions, it is necessary for $Q>0$. Following Ref.~ \cite{jawaharlal2024reduced}, it is known that when the entire trajectory is desensitized, the solution exhibits a lower dispersion in the cost value of the original OCP (i.e., in $m(t_f)$) in the presence of uncertainties, but desensitization results in an overall lower final cost value compared to the sensitive final cost value. In the context of fuel-optimal trajectories, the final mass of the spacecraft will inevitably decrease due to desensitization (i.e., $m(t_f)_{@Q= 0}>m(t_f)_{@Q>0}$). 

In this paper, we are investigating the impact of desensitization over specific time intervals of the trajectory and not along the entire trajectory. As an example, consider a single time interval defined as $t\in[t_1,t_2]$, where  $t_0 \leq t_1 < t_2 \leq t_f$. A time-triggered desensitization can be achieved by modifying the cost functional, given in Eq.~\eqref{eq:costfuncaug}, as,
\begin{equation} \label{eq:costfunctemporal}
     \underset{\delta \in [0,1] \& || \hat{\bm{u}}|| = 1}{\mathop{\text{minimize}}} ~J = -m(t_f) + \underbrace{\int_{t_0}^{t_f} \mu_1(t) \mu_2(t) Q \lambda^2_T(t) \ dt}_{\text{time-triggered desensitization penalty term}},
\end{equation}
where $\mu_1(t)$ and $\mu_2(t)$ are trigger factors that are defined as,
\begin{align} \label{eq:mu1mu2}
    \mu_1(t;t_1,\rho) = \frac{1}{2}  
    \left [ 1+\tanh \left (\frac{t-t_1}{\rho} \right ) \right ], \\
    \mu_2(t;t_2,\rho) = \frac{1}{2} \left [1+\tanh \left (\frac{t-t_2}{\rho} \right ) \right ],
\end{align}
where $t$ denotes the current time instant of the integration and $\rho$ is a smoothing parameter that controls the smooth transition of trigger factors from $0$ to $1$ and vice versa. The value for $\rho$ is set to $1.0 \times 10^{-5}$. Note that $\mu_1 \times \mu_2$ becomes equal to $1$, when $t_1 \leq t \leq t_2$, and becomes $0$ for all other values of $t$. 

\begin{figure}[ht!] 
    \centering
    \includegraphics[width=0.8\textwidth]{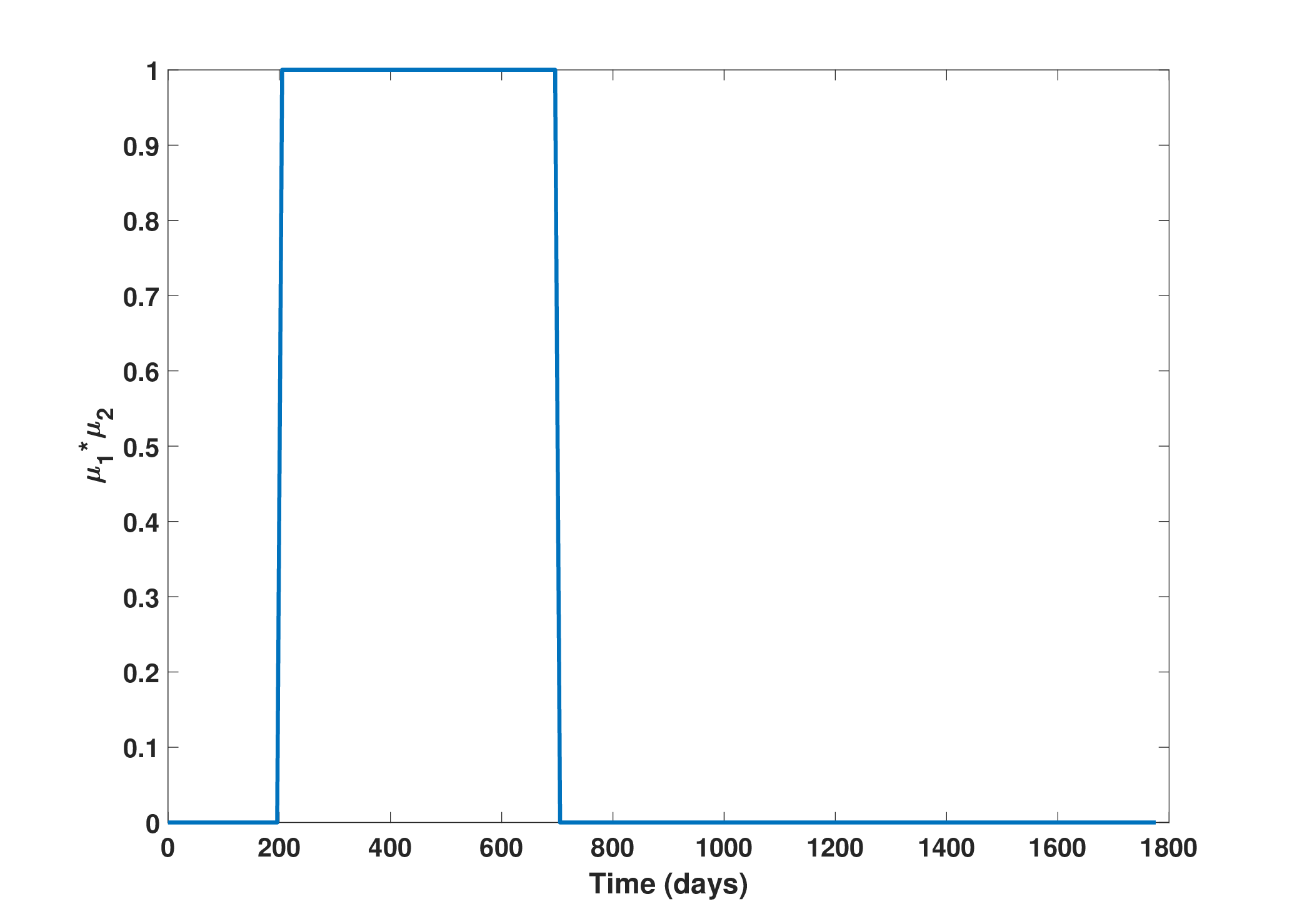}
    \caption{Time history of a representative $\mu_1 \times \mu_2$ with $t_1 = 200$ days and $t_2=700$ days.}
    \label{fig:mu1mu2}
\end{figure}

According to the time-triggered logic, the Lagrange cost in Eq.~\eqref{eq:costfunctemporal} becomes active only when $t_1 \leq t \leq t_2$. The idea of time-triggered activation is motivated through the Composite Smooth Control (CSC) method that allows for a smooth blending of multiple trigger factors to construct a signal \cite{taheri2020novel}. The CSC is a powerful method that has been used for solving complex OCPs such as enforcing forced-coast intervals prior and after gravity-assist maneuvers \cite{arya2021low}, enforcing duty-cycle constraints \cite{nurre2023duty}, eclipse constraints \cite{singh2021eclipse,taheri2021optimization,nurre2024end,sowell2024eclipse}, maximum payload indirect optimization of multimode propulsion systems \cite{arya2021composite,arya2021electric} and co-optimization of solar array size, trajectory and propulsion system operation modes using direct optimization methods \cite{saloglu2024co}. For example, Figure \eqref{fig:mu1mu2} shows the profile of $\mu_1 \times \mu_2$ when $t_0 = 0$, $t_1 = 200$ days, $t_2=700$ days and $t_f = 1776$ days. 

Following the procedure outlined in \cite{jawaharlal2024reduced}, the state dynamics, given Eq.~\eqref{eq:EOM}, and the thrust costate differential equation, given in Eq.~\eqref{eq:thrustcostateMEE}, are solved using a direct optimization method. We used a pseudo-spectral optimal control software, GPOPS-II \cite{patterson2014gpops}. We emphasize that the optimizer directly searches over the admissible set of controls. Since direct optimization over control inputs is performed, the following control constraints are also introduced to the problem formulation: $\delta \in [0,1]$ and $||\hat{\bm{u}}|| = 1$. The GPOPS-II setting is as follows: \textit{hp-LiuRao-Legendre} mesh method with a mesh tolerance of $1.0 \times 10^{-5}$ and IPOPT \cite{biegler2009large} is used as the non-linear programming (NLP) problem solver.

One of the advantages of using RDF is that desensitization can be achieved without using the costates of states and mass \cite{jawaharlal2024reduced}. Therefore, when GPOPS-II was used, the boundary conditions in Eq.~\eqref{eq:bcTPBVP} are modified as,
\begin{align} \label{eq:bcalt}
    \bm{\Psi} = \left[\ [\bm{x}(t_f) - \bm{x}_T]^\top, \lambda_T(t_f) \right]^\top = \bm{0}.
\end{align}

\section{Time-Triggered orbit-raising RDF formulation}\label{sec:orbraisrdf}

The time-triggered RDF is also applied to a classic orbit-raising problem, as discussed in \cite{bryson1975applied}. The state vector is denoted as $\bm{x} = [r,u,v]^{\top}$, where $r$ represents the radial coordinate, $u$ represents the radial velocity and $v$ represents the transverse velocity. The state dynamics for this problem are \cite{jawaharlal2024reduced}
\begin{align}\label{eq:OR_eq}
\dot{r} &= u, & \dot{u} &= \frac{v^2}{r} - \frac{\mu}{r^2} + \frac{T}{m_0 - |\dot{m}|t}\sin(\phi), & \dot{v} &= \frac{-uv}{r} + \frac{T}{m_0 - |\dot{m}|t} \cos(\phi), 
\end{align}
where $T$ denotes propulsive thrust, $m$ is the spacecraft mass, and $\phi$ is the thrust steering angle and denotes the scalar control input. The goal of this problem is to maximize the final radius and the time-triggered desensitizing cost functional is written as,
\begin{equation}
    \underset{\phi \in [0, 2 \pi]}{\mathop{\text{minimize}}}~J =  -r(t_f) + \underbrace{\int_{t_0}^{t_f}  \mu_1(t) \mu_2(t) Q \lambda_T^2 \ dt}_{\text{time-triggered desensitization penalty term}},
\end{equation}
where $\mu_1(t)$ and $\mu_2(t)$ have the same forms that are defined in Eq.~\eqref{eq:mu1mu2}.
Using the non-uniqueness property of the costates when  hybrid direct-indirect method is used, the thrust costate differential equation for this problem can be written as, \cite{jawaharlal2024reduced}
\begin{align} \label{eq:OR_CT}
  \dot{\lambda}_T = -\left(\frac{\sin(\phi)+\cos(\phi)}{m_0 + \Dot{m}(t-t_0)}\right).
\end{align}

Upon using the scaled values for states and parameters \cite{bryson1975applied}, the boundary conditions for this problem can be summarized as,
\begin{align} \label{eq:bcor}
    \bm{\Psi} = \left[\ [r(t_0) - r_0, \ u(t_0), \ v(t_0) - v_0, \ u(t_f), \ v(t_f)-\sqrt{\frac{1}{r(t_f)}}, \ \lambda_T(t_f) \right]^\top = \bm{0}.
\end{align}
\section{Results} \label{sec:results}
In this section, the results of time-triggered RDF for two classes of trajectory optimization problems are presented. First, we present the results for two fuel-optimal trajectories. Then, we present the solutions for orbit-raising types of maneuvers. 

\subsection{Earth-Comet 67P/Churyumov-Gerasimenko Fuel-Optimal Rendezvous Problem}
A fixed-time, fuel-optimal rendezvous maneuver from Earth to Comet 67P is considered \cite{praveen2021comparison}. The spacecraft's initial position vector is $\bm{r}(t_0) = [-10687809.15, -151602518.3, 8676.494013]^{\top}$ km and the initial velocity vector is $\bm{v}(t_0) = [29.22497601, -2.197707221, 0.000972199]^{\top}$ km/s. The target position vector  is $\bm{r}_T =  [-536251927.7, -126576922.3, 14541016.26]^{\top}$ \text{km} and the target velocity vector is $\bm{v}_T = [-6.858900316, -13.35248149, -0.453167946]^{\top}$ km/s. A spacecraft with an initial mass of $3000$ Kg is considered with the following propulsion system parameters: nominal initial thrust, $T_0=0.6$ N, constant specific impulse, $I_\text{sp}=3000$ s. 

The time-triggered RDF is applied to this problem with $t_1 = 0$. The value of $t_2$ is varied between $t_0$ and $t_f$. Based on the findings in Ref. \cite{jawaharlal2024reduced}, where the Earth-Comet 67P scenario demonstrated good desensitized results, the $Q$ value for this problem is fixed at 0.001.
Initially, the problem is solved for the nominal thrust value, $T_0 = 0.6$ N. To check for any dispersion, the problem is also solved with a $\pm \%5$ variation in initial thrust, i.e., for $T_0 = 0.63$ N and $T_0 = 0.57$ N. The dispersion level is evaluated by calculating the difference in the final mass ($m(t_f)$) when the thrust value is altered by five percent.

\begin{figure}[htbp!] 
    \centering
    \includegraphics[width=1\textwidth]{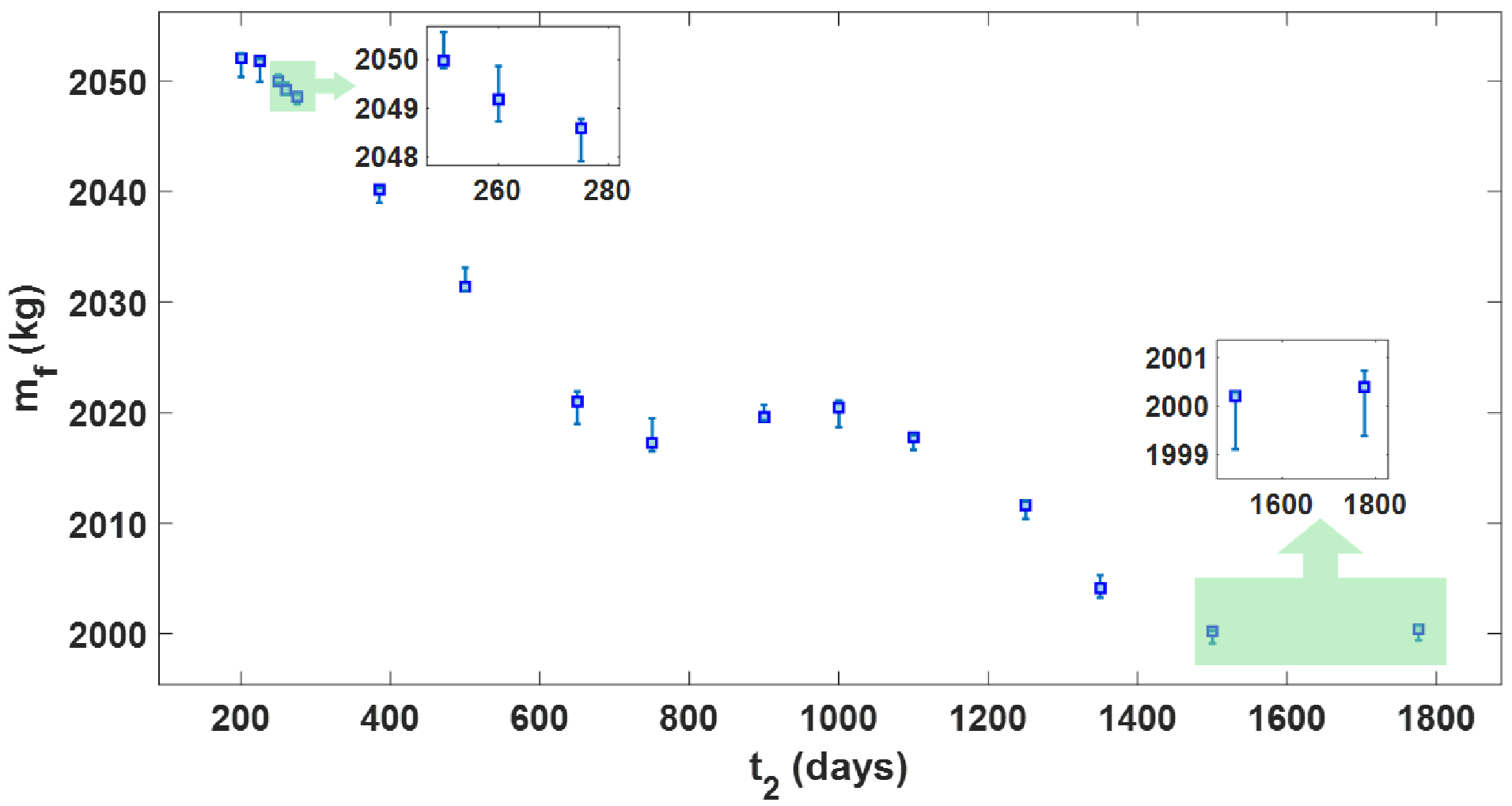}
    \caption{Fuel-optimal Earth-Comet67p problem: dispersion in $m_f$ vs. $t_2$ values.}
    \label{fig:EC_error}
\end{figure}

\begin{figure}[htbp!] 
    \centering
    \includegraphics[width=1.0\textwidth]{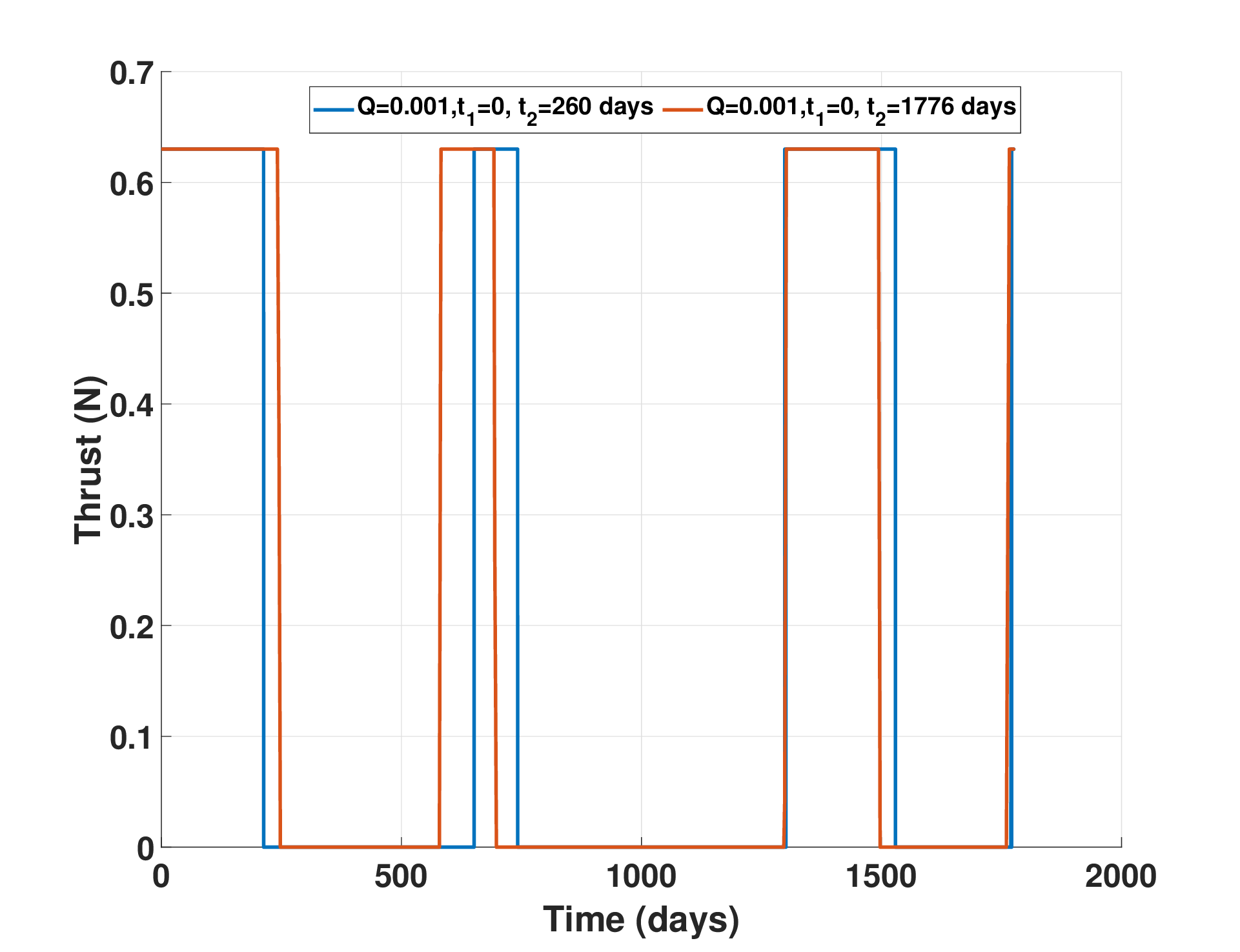}
    \caption{Earth-67P problem: desensitized thrust time histories for different $t_2$ values.}
    \label{fig:EC_thrust}
\end{figure}

\begin{figure}[htbp!] 
    \centering
    \includegraphics[width=1.0\textwidth]{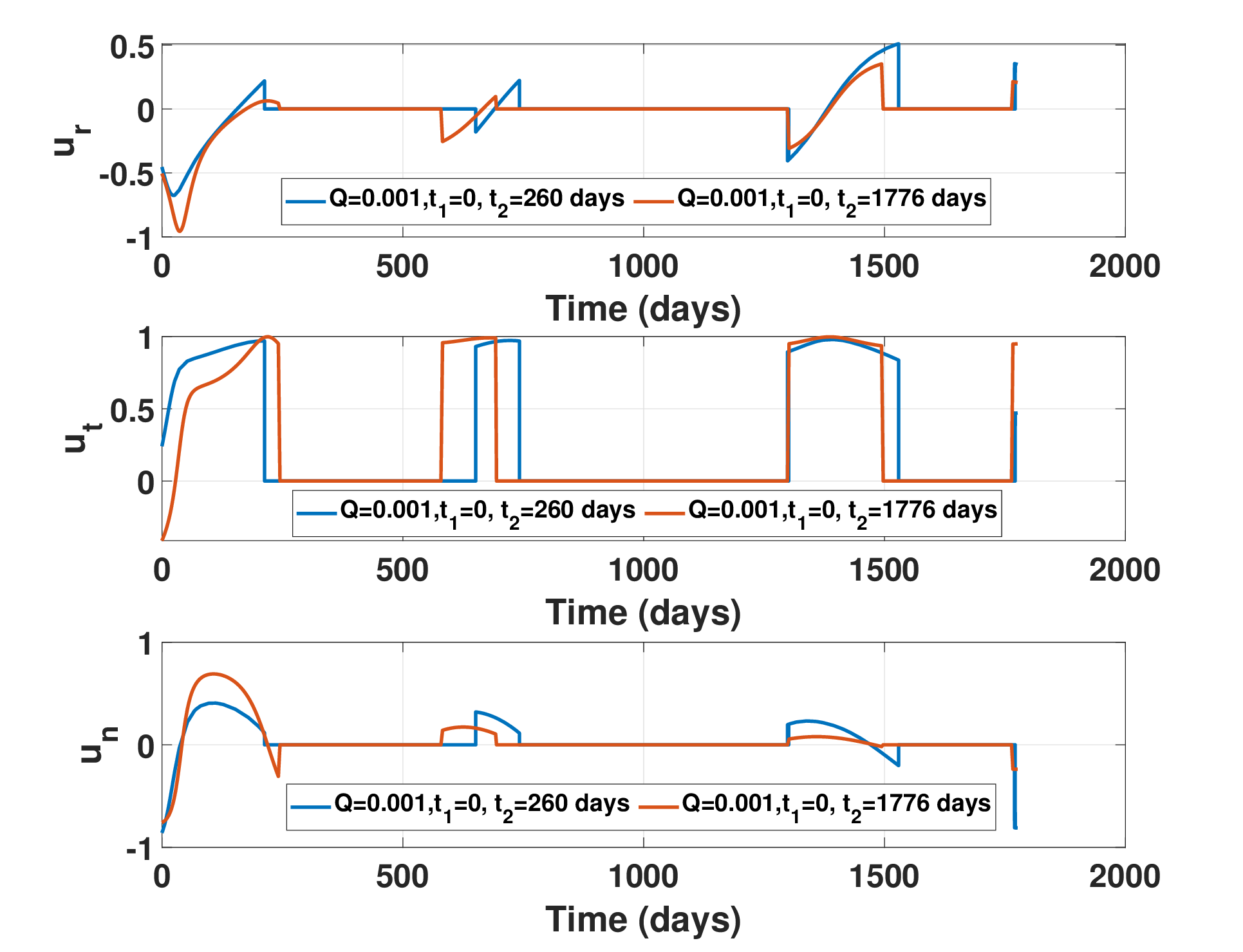}
    \caption{Earth-67P problem: desensitized thrust steering time histories for different $t_2$ values.}
    \label{fig:EC_controls}
\end{figure}

Figure \eqref{fig:EC_error} summarizes the dispersion in $m(t_f)$ as an error plot when $t_2$ is gradually increased till $t_2 = t_f$. The error plot shows both positive and negative dispersions, corresponding to increased and decreased thrust, respectively. Notably, the final mass does not decrease monotonically, with an increase observed near $t_2 \in [750,1000]$ days. Furthermore, certain values of $t_2$ yield a better final mass while maintaining similar dispersion levels compared to when the problem is desensitized for the entire time interval.
Table \ref{tab:EC_260} summarizes the dispersion levels for $t_2 = 260$ days and $t_2 = 1776$ days (total mission time), where $d$ represents the absolute difference in final mass ($m(t_f)$) when thrust is changed by $5$ percent. The dispersion level for $t_2 = 260$ days is comparable to that of the full mission desensitization, but with a significantly higher optimal final mass, approximately 49 kg greater.
Figures \eqref{fig:EC_thrust} and \eqref{fig:EC_controls} depict the evolution of thrust switching and components of the unit thrust direction vector, respectively, for both $t_2 = 260$ days and $t_2 = 1776$ days.

\begin{table*}[htbp!]
\caption{\textbf{Fuel-optimal Earth-67P problem: desensitization results for different $t_2$ values.}}
\label{tab:EC_260}
\centering 
\begin{tabular}{|c|c|c|c|c|c|}
\hline
$Q$ & $T_0$ (N) & $t_1$ (days) & $t_2$ (days) & $m(t_f)$ (kg) & $d$ (kg)    \\
\hline
0.0001 & 0.6 & 0 & 260& 2049.1819  &  -- \\
0.0001 & 0.63 &0 &260 & 2049.6331 & 1.1649   \\
0.0001 & 0.57 &0 &260 & 2049.8615 & 0.4937   \\
\hline
0.0001 & 0.6 & 0 & 1776 & 2000.5518 &  -- \\
0.0001 & 0.63 & 0 & 1776 & 1999.3869 & 0.4512  \\
0.0001 & 0.57 & 0 & 1776 & 2000.0582 & 0.6797  \\
\hline
\end{tabular}
\end{table*}

\subsection{Earth-Dionysus Fuel-Optimal Rendezvous Problem}
The spacecraft's initial state in Cartesian coordinates are (taken in the Sun-Centered Inertial frame), initial position: $\bm{r}(t_0) = [-3637871.081, 147099798.784, -2261.441]^{\top}$ km, initial velocity: $\bm{v}(t_0) = [-30.265097, -0.8486854, 0.0000505]^{\top}$ km/s. 
The target states in Cartesian coordinates are, target position: $\bm{r}_T = [-302452014.884, 316097179.632, 82872290.075]^{\top}$ km and target velocity: $\bm{v}_T = [-4.53347379984, -13.1103098008, 0.65616382602]^{\top}$ km/s. These values are converted into the set of MEEs for solving the OCPs. The mission and spacecraft parameters are as follows. Total time of flight: $3534$ days, initial spacecraft mass: $4000$ kg, nominal initial thrust: $0.32$ N and constant specific impulse ($I_\text{sp}$): $3000$ s. 

For the Earth-Dionysus problem, we implement the time-triggered RDF by setting $t_1 = 0$ and varying $t_2$ between $t_0$ and $t_f$. The $Q$ value for this problem is fixed as $0.0001$. Initially, the problem is solved for the nominal thrust value, $T_0 = 0.32$ N. To investigate the problem's sensitivity to thrust variations, we then solve the problem with a $\pm \%2$ thrust variation, using $T_0 = 0.3264$ N and $T_0 = 0.3136$ N. The dispersion is quantified by measuring the change in final mass ($m(t_f)$) resulting from this $2$ percent thrust variation.

\begin{figure}[htbp!] 
    \centering
    \includegraphics[width=1\textwidth]{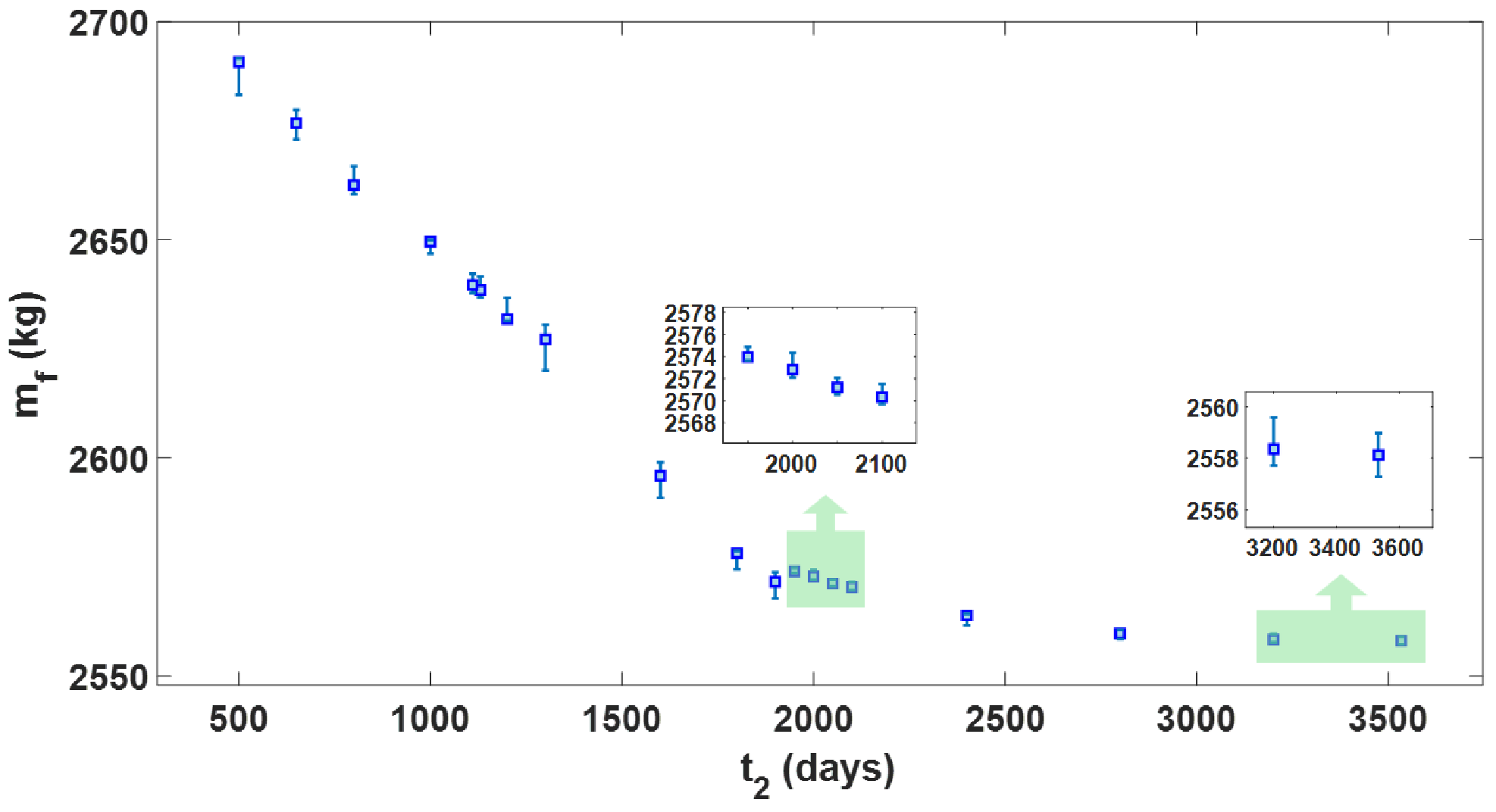}
    \caption{Earth-Dionysus problem: variation of $m_f$ for different $t_2$ values.}
    \label{fig:ED_error}
\end{figure}

Figure \eqref{fig:ED_error}, shows the dispersion of the final mass when the value of $t_2$ is gradually increased till $t_2 = t_f$. Notably, the relationship is non-monotonic, similar to the observations from the Earth-67P case. Additionally, certain values of $t_2$ yield a better final mass while maintaining similar dispersion levels compared to when the problem is desensitized for the entire time interval. 

In Table \ref{tab:ED_1950}, we compare dispersion levels for $t_2 = 1950$ days and $t_2 = 3534$ days (the full mission duration). Here, $d$ denotes the absolute difference in $m(t_f)$ when thrust increases by $2$ percent. The dispersion level, when the problem is desensitized for $1950$ days, is comparable to when the problem is desensitized for the entire mission time. However, the optimal final mass is higher when $t_2 = 1950$ days. 

\begin{figure}[ht!] 
    \centering
    \includegraphics[width=1.0\textwidth]{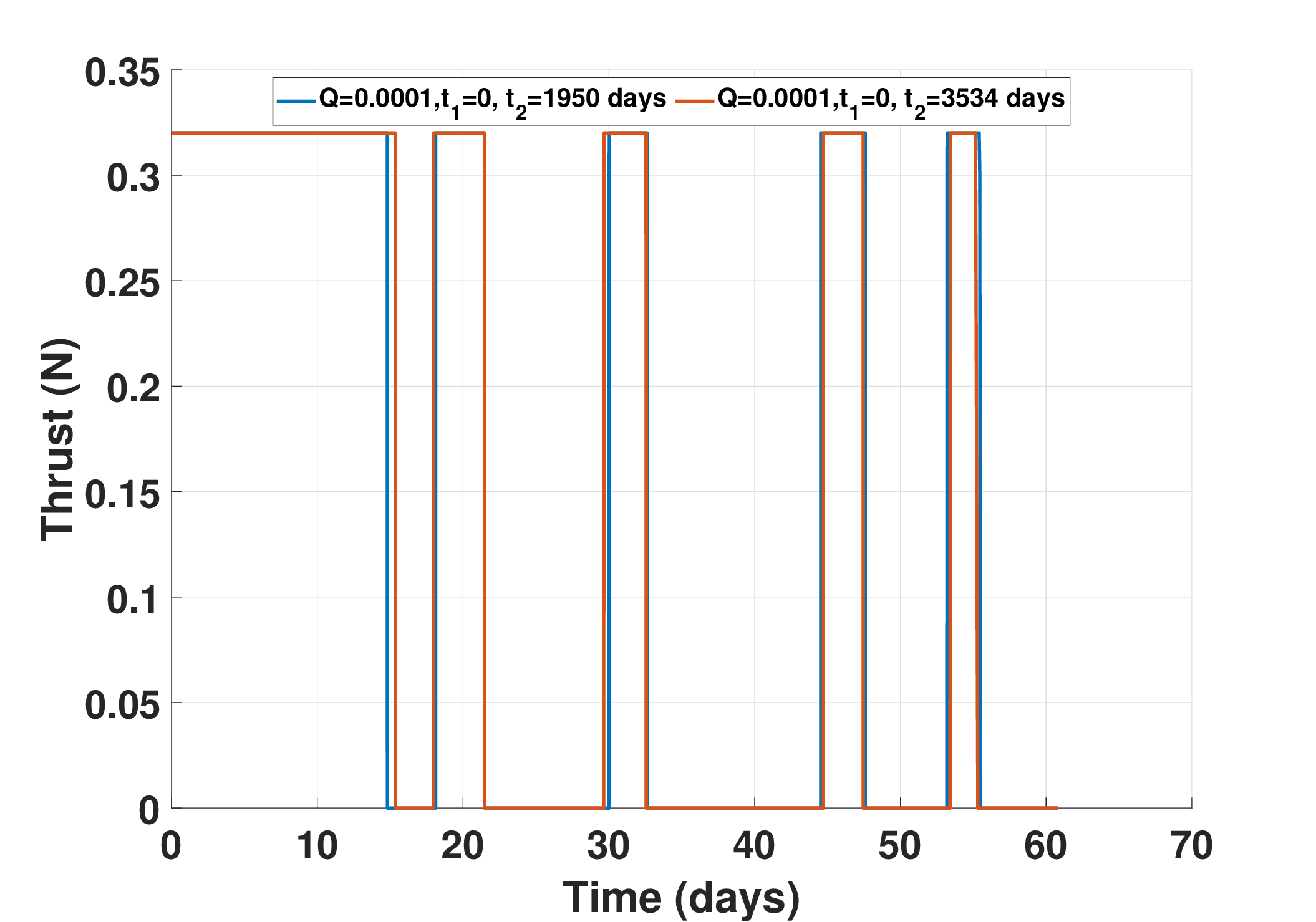}
    \caption{Earth-Dionysus problem: desensitized thrust vs. time for different $t_2$ values.}
    \label{fig:ED_T_switch}
\end{figure}

\begin{figure}[ht!] 
    \centering
    \includegraphics[width=0.9\textwidth]{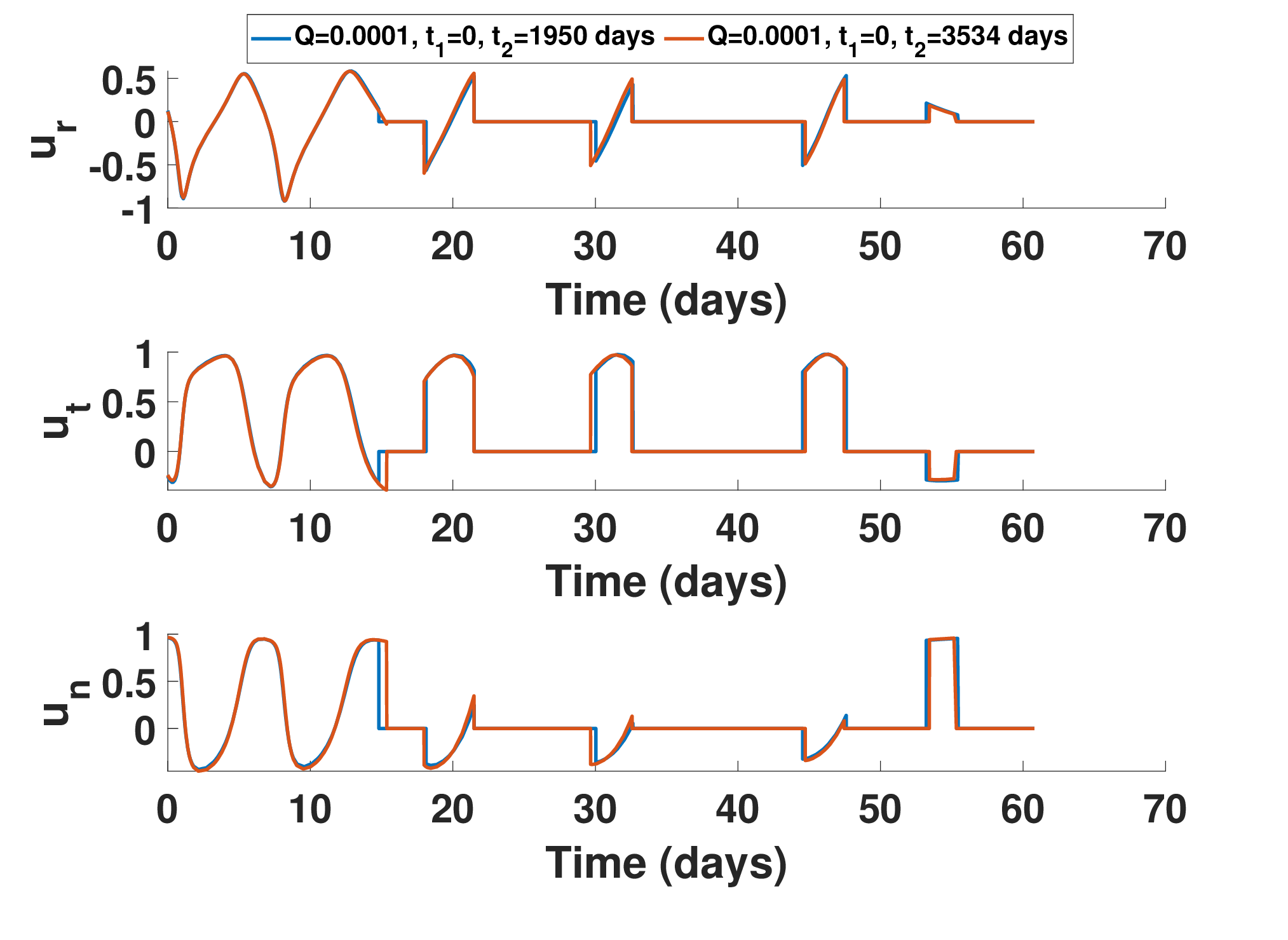}
    \caption{Earth-Dionysus problem: desensitized thrust steering vs. time for different $t_2$ values.}
    \label{fig:ED_controls}
\end{figure}

The plots for the evolution of thrust switching and components of the unit thrust direction vector are shown in Figures \eqref{fig:ED_T_switch} and \eqref{fig:ED_controls} for $t_2 = 1950$ days and $t_2 = 3534$ days, respectively.

\begin{table*}[htbp!]
\caption{\textbf{Earth-Dionysus problem: desensitization simulation results for different $t_2$ values.}}
\label{tab:ED_1950}
\centering 
\begin{tabular}{|c|c|c|c|c|c|}
\hline
$Q$ & $T_0$ (N) & $t_1$ (days) & $t_2$ (days) & $m(t_f)$ (kg) & $d$ (kg)    \\
\hline
0 & 0.32 & 0 & 1950 & 2573.9745  &  -- \\
0 & 0.3264 &0 & 1950 & 2574.3378 & 0.650   \\
0 & 0.3136 &0 & 1950 & 2573.1139 & 0.857
  \\
\hline
0.0001 & 0.32 & 0 & 3534 & 2558.1230 &  -- \\
0.0001 & 0.3264 & 0 & 3534 & 2558.9599 & 0.837  \\
0.0001 & 0.3136 & 0 & 3534 & 2557.2601 & 0.863  \\
\hline
\end{tabular}
\end{table*}

\begin{figure}[ht!] 
    \centering
    \includegraphics[width=1.1\textwidth]{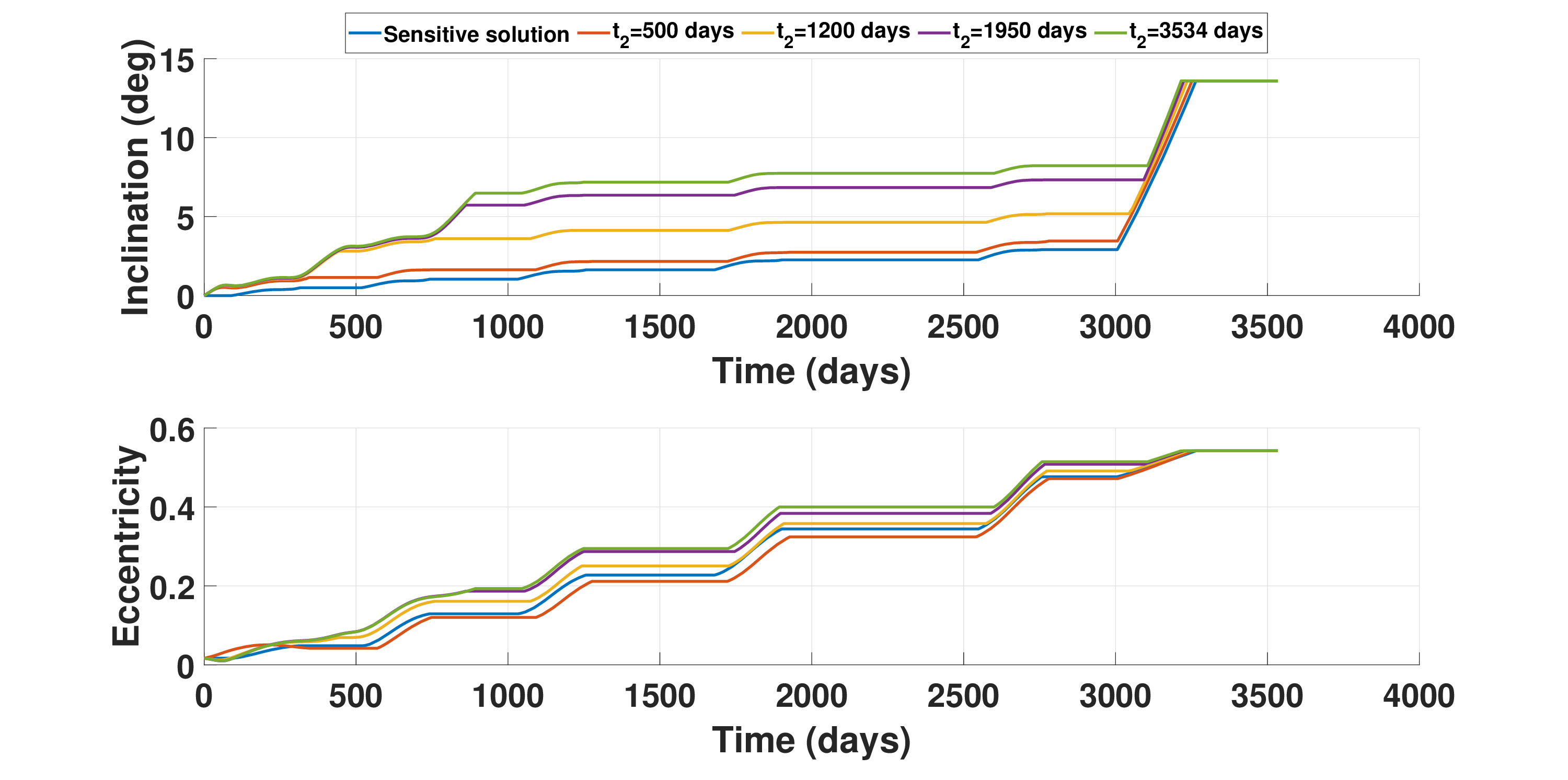}
    \caption{Earth-Dionysus problem: inclination and eccentricity vs. time for different $t_2$ values.}
    \label{fig:ED_OE}
\end{figure}

\newpage
\subsection{Orbit-Raising Problem}
The time-triggered desensitized results for the orbit-raising problem are summarized here. The scaled boundary conditions at the initial time $t_0$ are: $r_0=1$, $u_0=0$,
$v_0=1$, $m_0=1$, and $T_0 = 0.1405$. The problem is set with a fixed total flight time of 3.32 in scaled units. Additionally, the rate of change of mass is given as a fixed value, $\dot{m} = -0.0749$.

Similar to the previous cases, the value of $t_2$ is varied between $t_0$ and $t_f$ and $t_1 =0$. The $Q$ value for this problem is fixed as $0.0004$. First, this problem is solved for the nominal thrust value, $T_0=0.1405$ N. To check for dispersion in the final radius, this problem is also solved by increasing and decreasing the initial thrust to $T_0 = 0.1505 N$ and $T_0 = 0.1305 N$. The dispersion level is checked by calculating the difference in the mass at the final time ($m(t_f$)), when the thrust value is changed.

\begin{figure}[htbp!] 
    \centering
    \includegraphics[width=1\textwidth]{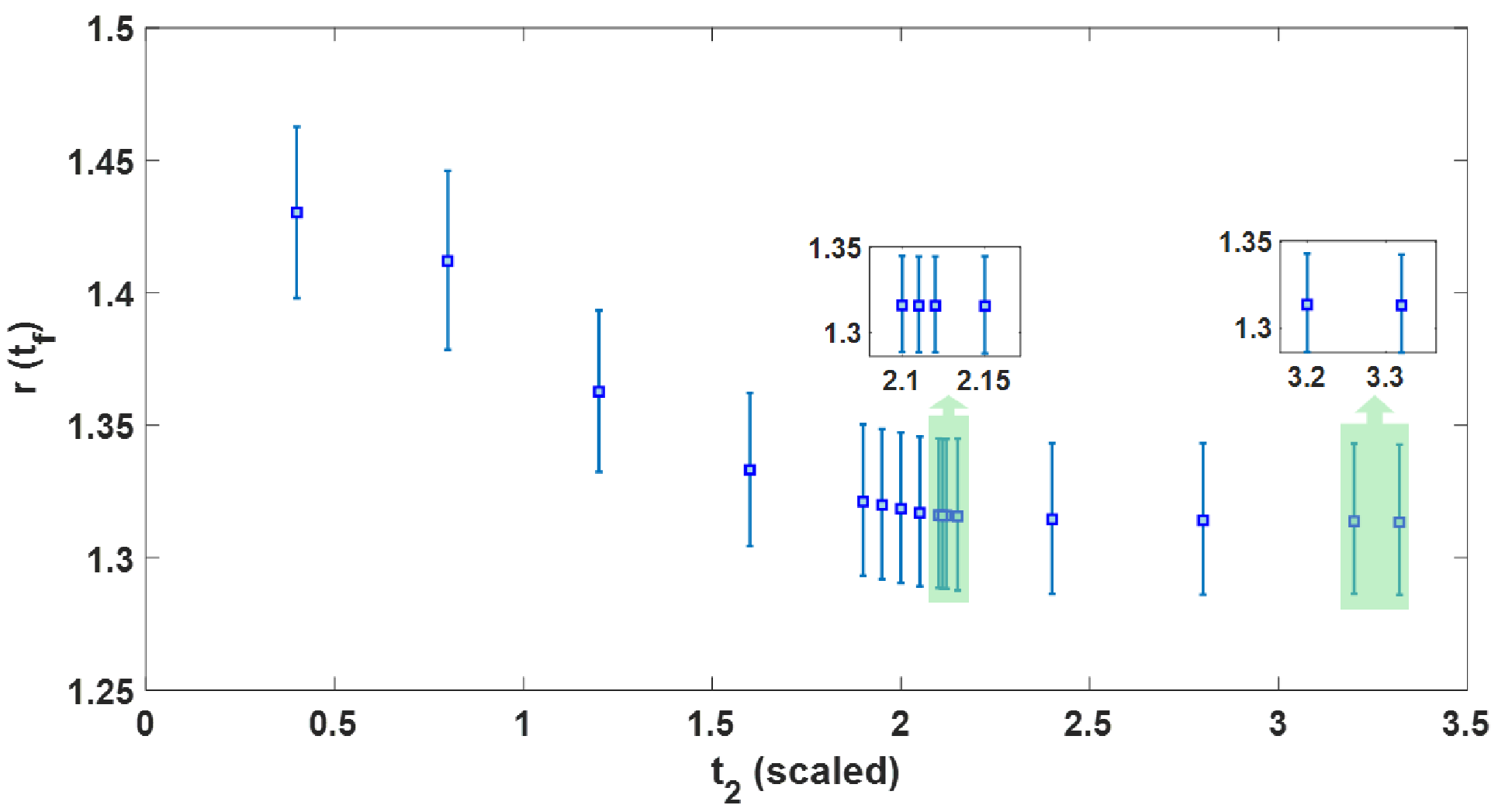}
    \caption{Orbit-raising problem: dispersion in $r (t_f)$ for different $t_2$ values.}
    \label{fig:OR_error}
\end{figure}

Figure \eqref{fig:OR_error}, shows the dispersion of the final radius when the value of $t_2$ is gradually increased till $t_2 = t_f$. There are some values of $t_2$ which have a better final radius while maintaining the same dispersion levels, when the problem is desensitized for the entire time interval. 

Table \ref{tab:OR_2p12} summarizes the dispersion levels for $t_2=2.12$ days and $t_2=3.32$ (total mission time). Here, $d$ is the absolute value of the difference in the radius at the final time ($r(t_f$)), when the value of thrust is changed.

\begin{table*}[htbp!]
\caption{\textbf{Orbit-raising problem: desensitization results for different $t_2$ values.}}
\label{tab:OR_2p12}
\centering 
\begin{tabular}{|c|c|c|c|c|c|}
\hline
$Q$ & $T_0$ (N) & $t_1$ (scaled)  & $t_2$ (scaled) & $r(t_f)$ & $d$     \\
\hline
0 & 0.1405 & 0 & 2.12  & 1.3158  &  -- \\
0 & 0.1505 &0 & 2.12 & 1.3447 & 0.029  \\
0 & 0.1305 &0 & 2.12 & 1.2883 & 0.028
  \\
\hline
0.0001 & 0.1405 & 0 & 3.32 & 1.3134 &  -- \\
0.0001 & 0.1505 & 0 & 3.32 & 1.3427 & 0.029  \\
0.0001 & 0.1305 & 0 & 3.32 & 1.2860 & 0.027 \\
\hline
\end{tabular}
\end{table*}

\begin{figure}[ht!] 
    \centering
    \includegraphics[width=1.0\textwidth]{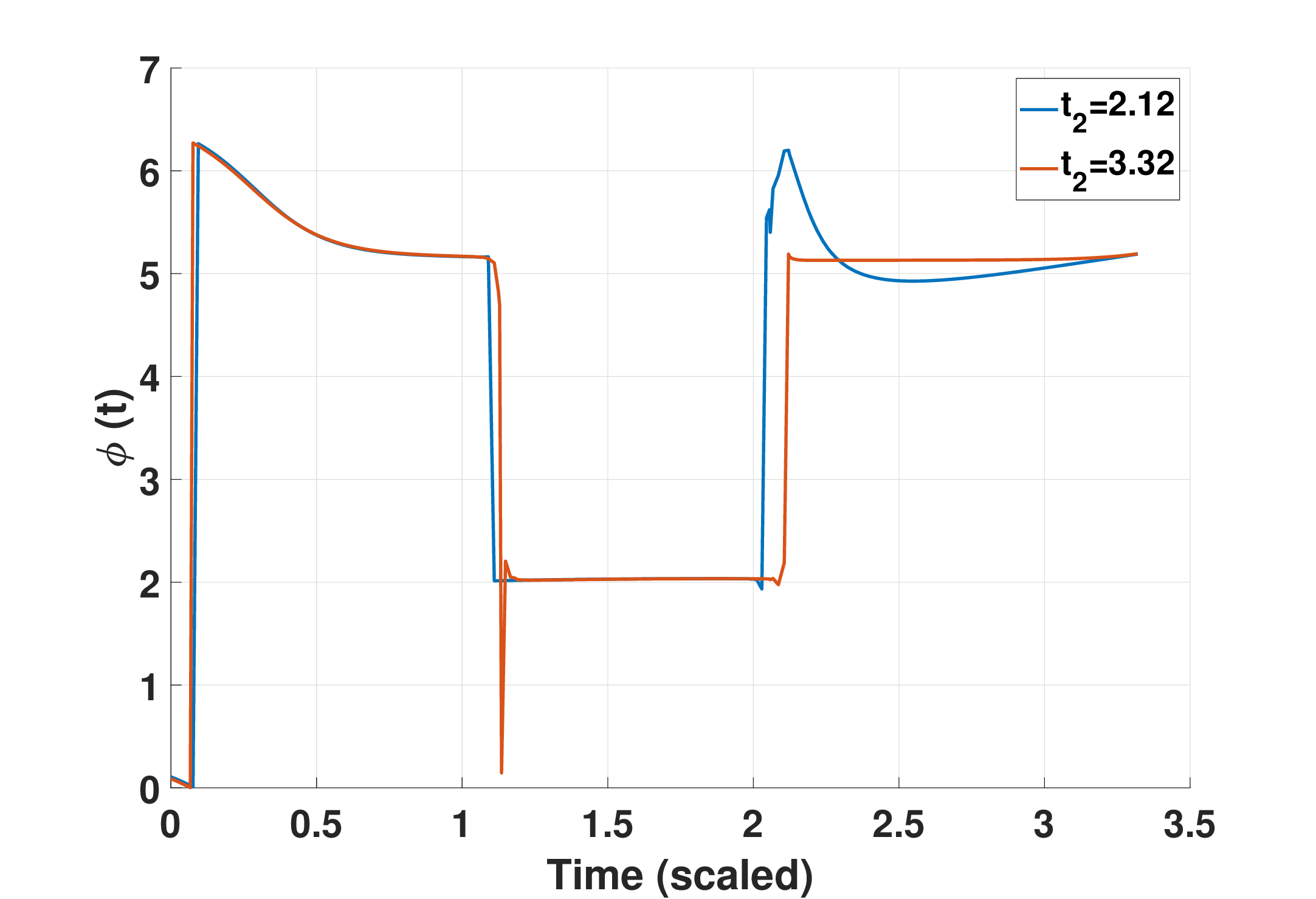}
    \caption{Orbit-raising problem: control vs. time for different values of $t_2$.}
    \label{fig:OR_control}
\end{figure}

\newpage
\section{Conclusion} \label{sec:conclusion}
We introduced a time-triggered Reduced Desensitization Formulation (RDF) for desensitizing optimal control problems. The time-triggered RDF selectively applies the desensitizing penalty term at specific intervals, thereby localizing desensitization to specific time intervals. The proposed time-triggered RDF formulation is then applied to the following trajectory optimization problems: 1) Earth-Comet 67P fixed-time fuel-optimal rendezvous maneuver, 2) Earth-Dionysus fixed-time fuel-optimal rendezvous maneuver, and 3) a classic orbit-raising problem. The numerical results suggest the following,
\begin{enumerate}
    \item Partial time interval desensitization can yield trajectories with improved final mass for fuel-optimal problems (or radius for orbit-raising problems) while maintaining dispersion levels comparable to full-trajectory desensitization.
    \item The relationship between desensitization interval and cost functional value is generally non-monotonic, suggesting the existence of optimal desensitization time windows.
    \item The time-triggered RDF approach offers flexibility in tailoring the desensitization strategy to specific mission intervals.
\end{enumerate}
Our future work focuses on extending this approach to multi-interval desensitization with other types of uncertainties, beyond thrust variations, could further enhance its practical utility in space mission design.

\bibliographystyle{AAS_publication}   
\bibliography{references.bib}   

\end{document}